\documentclass{article}
\usepackage{amsmath}
\usepackage{amssymb}
\usepackage{amsthm}
\newcommand{\ignore}[1]{\relax}

\newcommand{\C}{\mathbb C}
\newcommand{\R}{\mathbb R}
\newcommand{\Z}{\mathbb Z}

\newtheorem{thm}{Theorem}
\newtheorem{thmquote}{Theorem}

\newtheorem{lem}{Lemma}
\newtheorem{cor}{Corollary}

\theoremstyle{definition}
\newtheorem{defn}{Definition}
\theoremstyle{remark}

\newtheorem{rmk}{Remark}
\newcommand{\tor}{(\C^*)^2}
\newcommand{\rtor}{(\R^*)^2}

\newcommand{\conj}{\operatorname{conj}}
\newcommand{\dd}{\partial}
\newcommand{\am}{\mathcal{A}}

\newcommand{\Log}{\operatorname{Log}}
\newcommand{\grad}{\operatorname{grad}}
\newcommand{\Hess}{\operatorname{Hess}}

\newcommand{\MA}{\operatorname{M}}
\newcommand{\T}{\mathbf{T}}
\renewcommand{\setminus}{\smallsetminus}

\newcommand{\Area}{\operatorname{Area}}

\begin{document}

\title{Amoebas of maximal area.}
\author{Grigory Mikhalkin
\thanks{The first author is partially supported by the NSF}\\
\normalsize Department of Mathematics\\
\normalsize University of Utah\\
\normalsize Salt Lake City, UT 84112, USA\\
\normalsize mikhalkin@math.utah.edu
\and
Hans Rullg{\aa}rd\\
\normalsize Department of Mathematics\\
\normalsize Stockholm University\\
\normalsize S-10691 Stockholm, Sweden\\
\normalsize hansr@matematik.su.se}
%\author{Grigory Mikhalkin and Hans Rullg{\aa}rd}
%\thanks{The first author is partially supported by the NSF}
%\address{Department of Mathematics\\ University of Utah\\
%Salt Lake City, UT 84112, USA}
%\email{mikha@math.utah.edu}
%\address{Department of Mathematics\\ Stockholm University\\
%S-10691 Stockholm, Sweden}
%\email{hansr@matematik.su.se}
\maketitle

\begin{abstract}
To any algebraic curve $A$ in $\tor$
one may associate a closed infinite region $\am$
in $\R^2$ called {\em the amoeba} of $A$. The amoebas of
different curves of the same degree come in
different shapes and sizes.
All amoebas in $\rtor$ have finite area and,
furthermore, there is an upper bound on the
area in terms of the degree of the curve.

The subject of this paper is the curves in $\tor$
whose amoebas are of the maximal area.
We show that up to multiplication by a constant in $\tor$
such curves are defined over $\R$ and, furthermore, that
their real loci are isotopic to so-called {\em Harnack curves}.
\end{abstract}

\section{Introduction.}
Let $f:\C^2\to\C$ be a polynomial, $f(z_1,z_2)=\sum\limits_{j,k}
a_{jk}z_1^jz_2^k$. Its zero set in $\tor$ is a curve
$A=f^{-1}(0)\cap\tor$ (where $\C^*=\C\setminus\{0\}$).
Let $\Delta\subset\R^2$ be the Newton polygon of $f$,
i.e. the convex hull of $\{(j,k)\ |\ a_{jk}\neq0\}$.
Gelfand, Kapranov and Zelevinski introduced one more
object associated to $f$.
\begin{defn}[Gelfand, Kapranov, Zelevinski \cite{GKZ}]
\footnote{In this paper we restrict our attention to functions of two
variables. Amoebas are defined for functions of any number of variables.}
The amoeba $\am\subset\R^2$ of $f$ is $\Log(A)$,
where $\Log:\tor\to\R^2$, $(z_1,z_2)\mapsto(\log|z_1|,\log|z_2|)$.
\end{defn}
It was remarked in \cite{GKZ} that every component of $\R^2\setminus\am$
is open and convex in $\R^2$.
In particular, $\am$ is closed and its (Lebesgue) area
is well-defined.

Note that $\am$ is never bounded in $\R^2$, since $f^{-1}(0)$
must intersect the coordinate axes in $\C^2$.
However it was shown by Passare and Rullg{\aa}rd \cite{PaRu}
that the area of $\am$ is always finite.
Furthermore, it is bounded in terms of $\Delta$.
\begin{thmquote}[Passare, Rullg{\aa}rd \cite{PaRu}]
\begin{equation}
\label{areaineq}
\Area(\am)\le\pi^2\Area(\Delta).
\end{equation}
\end{thmquote}
The main result of this paper is the extremal property
of this inequality.
%\begin{thm}
%\label{max=>harnack}
%If $\Area(\am)=\pi^2\Area(\Delta)$ and $f$ is irreducible
%then there exist constants $a,b_1,b_2\in\C$ such that
%$g(z_1,z_2)=af(b_1z_1,b_2z_2)$ has real coefficients.
%Furthermore, the topological type of $((\R^*)^2,g^{-1}(0)\cap
%(\R^*)^2)$ depends only on $\Delta$, namely, the real locus of
%$g$ is a so-called Harnack curve, see Theorem \ref{ha}.
%\end{thm}

We say that a curve $A$ is defined over $\R$ if it is invariant under
the complex conjugation $\conj:\tor\to\tor$, $(z_1,z_2)\mapsto
(\bar{z}_1,\bar{z}_2)$. In this case
we may consider the real part of the curve $\R A=A\cap(\R^*)^2$
which is a real algebraic curve.
We say that a curve $A$ is real up to
multiplication by a constant if there exist constants
$b_1,b_2\in\C^*$ such that $(b_1,b_2)\times A\subset\tor$
is defined over $\R$.
The condition that $A$ is real up to multiplication by
a constant is equivalent to the condition that there exist
$a,b_1,b_2\in\C^*$ such that the polynomial $af(\frac{z_1}{b_1},
\frac{z_2}{b_2})$ has real coefficients.
In this case we may also consider the real part
$\R A=\{(x_1,x_2)\in\rtor\ |\ af(\frac{x_1}{b_1},\frac{x_2}{b_2})=0\}$.
We say that a map is at most 2-1 if the inverse image
of any point in the target consists of at most 2 points.
The main result of this paper is the following theorem.
\begin{thm}
\label{main}
Suppose that $\Area(\Delta)>0$.
Then the following conditions are equivalent.
\begin{enumerate}
\item $\Area(\am)=\pi^2\Area(\Delta).$
\item The map $\Log|_A: A\to\R^2$ is at most 2-1
and $A$ is real up to multiplication by a constant.
\item The curve $A$ is real up to multiplication
by a constant and its real part $\R A$ is a (possibly
singular) Harnack curve
(see Definitions \ref{ha} and \ref{sha}) for the Newton polygon $\Delta$.
\end{enumerate}
Furthermore, these conditions imply that the non-singular locus
of $\R A$ coincides with $A\cap\Log^{-1}(\dd\am)$.
\end{thm}

\begin{cor}
The inequality \eqref{areaineq} is sharp for any Newton polygon $\Delta$.
\end{cor}
The corollary follows from Theorem \ref{main} and the Harnack-Itenberg-Viro
Theorem (see section \ref{harnackcurves}) on existence of Harnack curves.

\begin{rmk}
A curve that is real up to multiplication by a constant
may have more than one real part (other real parts may come
as a result of multiplication by different constants).
For instance, if $f$ is a real polynomial which contains only
even powers of $z_2$ then the pullback of $f$ under
$(z_1,z_2)\mapsto (z_1,iz_2)$ is a real polynomial with
a different real part.

The theorem implies that a Harnack
curve is real up to multiplication
by a constant in a unique way.
Indeed, the choice of the real part is determined
by the identity $\R A=A\cap\Log^{-1}(\dd\am)$.
\end{rmk}

%Recall (see e.g. \cite{GKZ}) that the
%polygon $\Delta$ determines a toric surface $\R T_\Delta\supset(\R^*)^2$.
%Let $\R Z=\R T_\Delta\supset(\R^*)^2$, it is a union of $n$ lines,
%$n$ ``coordinate axes" of $\R T_\Delta$, if $\Delta$ is an $n$-gon.
%The statement of Theorem \ref{max=>harnack} can be strengthened.
%
%\begin{add}
%If ...the topological type of the triad $(\R T_\Delta;\R Z,\R A)$
%\end{add}
%In the remaining part of the section we describe Harnack curves.

%\begin{exa}[Harnack curves, \cite{Ha}, \cite{IV}]
%\label{ha}
%\mnote{Description of Harnack curves}
%\section{Curves in $\rtor$ and in real toric surfaces.}
\section{Harnack curves in $\rtor$.}
\label{harnackcurves}
Let us fix a convex polygon $\Delta\subset\R^2$ whose
vertices have integer coordinates.
Consider all
possible real polynomials $f$ whose Newton polygon is $\Delta$.
The same polynomial $f$ may be viewed both as a function
$\tor\to\C$ and as a function $\rtor\to\R$.

Let $\R A$ be the zero set of $f$ in $\rtor$.
Equivalently, $\R A$ is a real part of the zero set $A$ of
$f$ in $\tor$.
For a generic choice of coefficients of $f$ the curve
$\R A$ is smooth.
However the topology of $((\R^*)^2,\R A)$ is different
for different choices of coefficients of $f$.
In particular, the number of components of $\R A$ may be
different. Also the mutual position of the components may
be different.

We may compactify the above setup.
Recall (see e.g. \cite{GKZ}) that the polygon $\Delta$
determines a toric surface $\C T_\Delta\supset\tor$. We denote the
real part of $\C T_\Delta$ with $\R T_\Delta\supset\rtor$.
The surface $\C T_\Delta$ is a compactification of $\tor$.
Furthermore, the complement $\C T_\Delta\setminus\tor$ is a union of $n$
(non-disjoint) lines, where $n$ is the number of sides of $\Delta$.
Similarly, $\R T_\Delta\setminus\rtor$ is a union of $n$ real lines
$l_1,\dots,l_n$.
These lines are called the {\em axes} of $\R T_\Delta$.
We assume that the indexing of $l_k$ is consistent with
the natural cyclic order on the sides of $\Delta$.

The closure $\bar{A}$ of $A\subset\tor\subset\C T_\Delta$
in $\C T_\Delta$ is a compact curve whose real part is
$\R\bar{A}\supset\R A$.
The topology of the triad $(\R T_\Delta;\R\bar{A},l_1,\cup\dots\cup l_n)$
carries all topological information on arrangement of $\R A$ in $\rtor$.

The upper bound on the number of components of $\R\bar{A}\subset
\R T_\Delta$ is provided
by Harnack's inequality \cite{Ha}. This number is never greater
than one plus the genus of $A$. Recall that by \cite{Kh}
the genus of $A$ is equal to the number of lattice points in the
interior of $\Delta$. We denote this number with $g$.

To deduce the upper bound on the number of components
of $\R A\subset\rtor$ we recall that $\R A=\R\bar{A}\setminus(l_1\cup\dots
\cup l_n)$, where $l_k$ corresponds to a side $\delta_k$
of $\Delta$.
Let $d_k$ be the integer length of $\delta_k$, i.e. the number
of lattice points inside $\delta_k$ plus one. Note that this
length is an $SL(2,\Z)$-invariant. The curve $\R\bar{A}$ and the axis
$l_k$ intersect in no more that $d_k$ points, since $d_k$ is
the intersection number of their complexifications. Therefore,
$\R A$ has no more than $g+\sum\limits_{k=1}^n d_k$ components.

\begin{defn}[Harnack curves, cf. \cite{Mi}]
\label{ha}
A non-singular curve $\R A\subset\rtor$ with the Newton polygon $\Delta$
is called a {\em Harnack curve} if all the following conditions hold.
\begin{itemize}
\item The number of components of $\R\bar{A}$ is equal to $g+1$
(where $g$ is the number of lattice points in the interior of $\Delta$).
%\item The number of components of $\R A$ is equal to
%$g+\sum\limits_{k=1}^n d_k$.
\item All components of $\R\bar{A}$ but one do not intersect
$l_1\cup\dots\cup l_n$.
\item A component $C$ of $\R\bar{A}$ can be divided into $n$
consecutive (with respect to the cyclic order on $C$) arcs
$\alpha_1,\dots,\alpha_n$ so that for each $k$ the intersections
$\alpha_k\cap l_k$ consists of $d_k$ points, while $\alpha_k\cap
l_j=\emptyset$, $j\neq k$.
\end{itemize}
\end{defn}
Note that the first two conditions imply that
the number of components of a Harnack curve $\R A$
is equal to $g+\sum\limits_{k=1}^n d_k$.

\begin{thmquote}[Mikhalkin \cite{Mi}]
For each Newton polygon $\Delta$ the topological type of the triad
$(\R T_\Delta;\R\bar{A},l_1\cup\dots\cup l_n)$ is unique
if $\R A$ is a Harnack curve.
\end{thmquote}

Note that the above theorem implies that the topological
type of the pair $(\rtor, \R A)$ is also unique for each $\Delta$.

\begin{thmquote}[Harnack, Itenberg, Viro, \cite{Ha}, \cite{IV}, \cite{Mi}]
Harnack curves exist for any Newton polygon $\Delta$.
\end{thmquote}

Harnack \cite{Ha} proved this theorem for
plane projective curves of arbitrary degree $d$.
In our language this corresponds to the case when
$\Delta$ is a triangle whose vertices are $(0,0)$, $(d,0)$,
$(0,d)$. Harnack's example was generalized to arbitrary Newton
polyhedra $\Delta$ with the help of Viro's patchworking
described in \cite{IV}, see Corollary A4 in \cite{Mi}.
The Harnack curves are a special case of the so-called T-curves,
see \cite{IV}.

We refer to \cite{IV} and \cite{Mi} for illustrations
of Harnack curves.
%In the notations of \cite{IV} a Harnack curve is a T-curve such that
%the sign of a node $(a,b)$ is $(-1)^{(a+1)(b+1)}$.

Recall that a point $p\in\R A\subset\rtor$ is called an {ordinary
real isolated double point} of $\R A$
(or an $\operatorname{A_1^+}$-point, see \cite{AVGZ})
if there exist local coordinates $x_1,x_2$ at $p\subset\rtor$
such that $A$ is locally defined by equation $x_1^2+x_2^2=0$.
%Alternatively, $p$ is the intersection point of a pair of
%transverse conjugate non-singular branches of $A$ ($x_2=ix_1$
%and $x_2=-ix_1$ in the local coordinate system above).

\begin{defn}[Singular Harnack curves]
\label{sha}
A singular curve $\R A\subset\rtor$ with the Newton polygon $\Delta$
is called a singular {\em Harnack curve} if
\begin{itemize}
\item the only singular points of $\R A$ are
$A_1^+$-points (ordinary real isolated double points);
\item the result of replacing of the singular points of $\R A$
with small ovals (which corresponds to replacing with the locus
$x_1^2+x_2^2=\epsilon, \epsilon>0$ in the local coordinates) gives
a Harnack curve for $\Delta$.
%a curve $\mathcal C\subset\rtor$ isotopic to a non-singular
%Harnack curve for $\Delta$.
\end{itemize}
\end{defn}
%\mnote{What's the difference in saying that a curve is isotopic to a
%Harnack curve and saying that it actually is Harnack?}
In other words, a singular Harnack curve is the result of
contraction to points of some ovals of a non-singular
Harnack curve.

\section{Monge-Amp\`ere measure on $\am$.}
\label{hans}

In the next section we prove the equivalence of conditions 1 and 2 in
the main theorem. The proof is an extension of the proof of the inequality
\eqref{areaineq} given in \cite{PaRu}.
We recapture in this section the main points in this proof.
The idea is to construct a measure on the amoeba  $\am$,
whose total mass is related to $\Delta$ and which can be
computed explicitly in terms of the hypersurface $A$.
This measure will be obtained as the real Monge-Amp\`ere
measure of a certain convex function associated to $f$.

We indicate briefly the definition of the real Monge-Amp\`ere
operator. Details may be found in \cite{RaTa}. Suppose $u$ is a smooth
convex function defined in $\R^n$. Then $\grad u$ defines a mapping from
$\R^n$ to $\R^n$. The Monge-Amp\`ere measure $\MA u$ of $u$ is defined by
$\MA u(E) = \lambda(\grad u(E))$ for any Borel set $E$, where $\lambda$
denotes Lebesgue measure on $\R^n$. That this is actually a measure requires
a proof, since $\grad u$ is in general not 1-to-1. If $u$ is convex but
not necessarily smooth, $\grad u$ can still be defined as a multifunction,
and the Monge-Amp\`ere measure of $u$ is defined as in the smooth case.
For smooth functions the Monge-Amp\`ere measure
is given by the determinant of the Hessian matrix,
$$\mu=|\Hess(u)|\lambda,$$
where $\lambda$ is the Lebesgue measure.

\paragraph{}Suppose now that $f$ is a given polynomial in two variables
and define
\begin{equation*}
N_f(x) = \frac{1}{(2\pi i)^2} \int_{\Log^{-1}(x)}
\frac{\log|f(z)|\, dz_1\, dz_2}{z_1 z_2}.
\end{equation*}
This is a real-valued function defined in $\R^2$, which is convex because
$\log|f(z)|$ is plurisubharmonic. Define $\mu$ to be the Monge-Amp\`ere
measure of $N_f$.

\begin{lem}\label{lem1}
The measure $\mu$ has its support in $\am$ and its total mass is
equal to the area of $\Delta$.
\end{lem}

\begin{proof} It is not difficult to show that $N_f$ is affine linear in
each connected component of $\R^2 \smallsetminus \am$ and that the
gradient image $\grad N_f (\R^2)$ is equal to $\Delta$ minus some of its
boundary points. This readily implies the statement. For details we refer
to \cite{PaRu}.
\end{proof}

\paragraph{}Let $F$ denote the set of critical values of the mapping
$\Log: A \to \R^2$. Pick a point $x_0 \in \am \smallsetminus F$
and functions $\phi_{j}, \psi_{j}$ defined in a neighborhood $V$ of $x_{0}$,
where $j$ ranges from 1 to $n$ and $n$ is the cardinality of
$\Log^{-1}(x_0) \cap A$,
such that
$A \cap \Log^{-1}(V) = \cup_{j = 1}^n
\{(\exp(x_{1} + i\phi_{j}(x)), \exp(x_{2} + i\psi_{j}(x)));x =
(x_{1},x_{2}) \in V\}$.
The main step in the proof of the inequality is the following computation.

\begin{lem}\label{lem2}
       With notations as above we have
\begin{equation}\label{hessian}
\Hess(N_{f}) = \frac{1}{2\pi} \sum_{j=1}^n
\pm \begin{pmatrix}\partial \psi_{j}/\partial x_{1} &
\partial \psi_{j}/\partial x_{2}\\ -\partial \phi_{j}/\partial x_{1} &
-\partial \phi_{j}/\partial x_{2} \end{pmatrix}.
\end{equation}
The signs depend on the signs of the intersection numbers between
$\Log^{-1}(x_0)$ and $A$. Each term in the sum is a symmetric, positive
definite matrix with determinant equal to 1.
\end{lem}

\paragraph{} For the proof we refer to \cite{PaRu}. We remark that the
fact that the matrices are symmetric with determinant equal to 1
follows immediately when we know that $A$ is a complex analytic curve.
The two last lemmas immediately imply the inequality \eqref{areaineq}
via the following corollary.

\begin{cor}\label{cor1}
If $\lambda$ denotes Lebesgue measure in $\R^2$, then
$\mu \geq (\lambda/\pi^2)|_{\am}$. Hence the area of $\am$ is not
greater than $\pi^2$ times the area of $\Delta$.
\end{cor}

\begin{proof}
It is not difficult to show that for $2 \times 2$  symmetric,
positive definite matrices $M_1, M_2$ the inequality
\begin{equation}\label{detineq}
\sqrt{\det(M_1 + M_2)} \geq \sqrt{\det M_1} + \sqrt{\det M_2}
\end{equation}
holds, with
equality precisely if $M_1$ and $M_2$ are real multiples of each other.
Applying this to the sum \eqref{hessian} and using the fact that it
contains at least two terms for all $x_0 \in \am \smallsetminus F$, the first
statement follows. Combining this with Lemma \ref{lem1} yields the second
part.
\end{proof}

\begin{rmk} The inequality used in the previous proof follows as a special
case of an inequality for positive definite matrices of arbitrary size,
analogous to the Alexandrov-Fenchel inequality for mixed volumes. The general
inequality can be found in \cite{Al}.
\end{rmk}

\section{Proof of Theorem 1: conditions 1 and 2 are equivalent.}

\paragraph{} We are now ready to prove the equivalence of conditions 1 and 2.
Note that by Corollary \ref{cor1},
$\Area(\am) = \pi^2 \Area(\Delta)$ if and only if
$\mu = (\lambda/\pi^2)|_{\am}$.

\subsection{Implication $1 \implies 2$.}
Suppose that $\mu = (\lambda/\pi^2)|_{\am}$. We first show that $f$ is
irreducible.

\begin{lem}\label{irreducible}
If $\mu = (\lambda/\pi^2)|_{\am}$, then $f$ is irreducible.
\end{lem}

\begin{proof} Let $K, L$ be compact convex subsets of $\R^2$. From the
monotonicity properties of mixed volumes it follows that
$\Area(K + L) \geq \Area(K) + \Area(L)$ with strict inequality holding
unless one of $K, L$ is a point or $K$ and $L$ are two parallel segments.
Assume now that we have a non-trivial factorization $f = gh$ and
let $\Delta_g, \Delta_h$ denote the Newton polytopes and $\am_g, \am_h$
the amoebas of $g$ and $h$ respectively. From Lemma \ref{lem1} it follows
that $\Area(\am) = \pi^2 \Area(\Delta)$. On the other hand, since
$\am = \am_g \cup \am_h$ and $\Delta = \Delta_g + \Delta_h$, it follows
from Corollary \ref{cor1} that
\begin{equation*}
\Area(\am) \leq \Area(\am_g) + \Area(\am_h) \leq \pi^2(\Area(\Delta_g) +
\Area(\Delta_h)) < \pi^2 \Area(\Delta).
\end{equation*}
This is a contradiction. \end{proof}

      From \eqref{detineq} it follows that for equality to
      hold in Corollary \ref{cor1} it is necessary that $\Log^{-1}(x)$
intersects $A$ in at most two points for all $x \not \in F$.
Hence the sum \eqref{hessian} contains two terms with
opposite signs. For equality to hold in \eqref{detineq} applied to the
sum \eqref{hessian}
it is necessary that $\grad \phi_1 = - \grad \phi_2$ and
$\grad \psi_1 = - \grad \psi_2$. After a multiplication of each
coordinate by a constant
we may
assume that $\phi_1 = - \phi_2, \psi_1 = - \psi_2$ in a neighborhood
of a given point in $\am \smallsetminus F$. (The existence of such points is
guaranteed by the assumption that $\Area(\Delta)$ and hence $\Area(\am)$ is
       positive.) But then $f(z)$ and
$\overline{f(\bar{z})}$ have a common factor, and hence coincide up to
a multiplicative constant since they are irreducible. Multiplying $f$ by
a suitable constant, we obtain a polynomial with real coefficients.

To complete the proof we must show
that $\Log^{-1}(x_{0})$ intersects $A$ in at most two points for all
$x_0 \in F$.
Note that $\Log^{-1}(x_0)\cap A$ cannot contain more
than 2 isolated points. Indeed, a small neighborhood in $A$ of an
isolated point in $\Log^{-1}(x_0)\cap A$ is mapped by $\Log$ either onto
a neighborhood of $x_{0}$, or in a 2-to-1 fashion onto a half-disk
with $x_{0}$ on its boundary. In any case, the presence of more than
2 isolated points would imply that $\Log^{-1}(x)\cap A$ contains more
than two points for some $x\notin F$, which is a contradiction.
%Otherwise there would exist $x\notin F$
%with $\Log^{-1}(x)\cap A$ consisting of more than two points
%since $A$ is connected.
%Let $x_0 \in F$ be given.
%If $\Log^{-1}(x_0) \cap A$ is
%discrete this is an easy argument similar to the considerations in
%the proof that $2 \implies 3$.

%Assume first
%that $\Log^{-1}(x_{0}) \cap A$ is a discrete set containing at least three
%points $p_1, p_2, p_3$. Let $V$ be a small neighborhood of $x_0$ and let
%$U_j$ be the connected component of $\Log^{-1}(V) \cap A$ containing $p_j$.
%Here we may assume that $U_1, U_2, U_3$ are mutually disjoint. Write
%$V_j = \Log(U_j)$. Now consider the following cases. If all $V_j$ contain
%a neighborhood of $x_0$, then $\Log^{-1}(x) \cap A$ has at least three
%points for all $x$ in this neighborhood, a contradiction. If one of the
%sets, say $V_1$, contains a neighborhood of $x_0$ and one, say $V_2$ does
%not, then $\Log^{-1}(x) \cap U_2$ contains at least two points for generic
%$x \in V_2$ and we once again obtain a contradiction.
%\mnote{Is it necessary to consider the above
%cases? It seems that by the rest of the proof $x_0\in F$ cannot
%have more than one inverse image.}
%Finally, if none of
%the sets $V_j$ contains a neighborhood of $x_0$, then by Lemma 1 in
%\cite{Mi} each  $V_j$ is concave in a neighborhood of $x_0$. Hence the
%intersection of two of these sets has non-empty interior. For generic $x$
%in this intersection, $\Log^{-1}(x) \cap A$ has at least 4 points, once
%again a contradiction.

If $\Log^{-1}(x_0) \cap f^{-1}(0)$ contains a curve $\gamma$ we
consider two different cases. If $\gamma$ is of the form
$\Log^{-1}(x_0) \cap \{z_1^jz_2^k = c\}$ for some $(j,k) \in \Z^2$ and
$c \in \C$, then $f$ contains the factor $z_1^jz_2^k - c$, which is impossible
by Lemma \ref{irreducible}. Otherwise, $t \gamma := \{(t_1 z_1, t_2 z_2);
(z_1,z_2) \in \gamma\}$ intersects $\gamma$ for all $t$ in an open set in
the real torus $\T^2$. By Theorem 5 in \cite{PaRu} (cf. the proof of
Lemma \ref{lem3})
this implies that $\mu$ has a point mass at $x_0$,
contradicting the assumptions. Hence we have shown that
$\Log: A \to \R^2$ is at most 2-to-1.

\subsection{Implication $2 \implies 1$.}
Conversely, assume that $\Log:A \to \R^2$ is at most 2-to-1 and that
$f$ has real coefficients. Since $\am$ and $\mu$ are invariant under the
changes of variables permitted in the theorem, this is no loss of generality.
Then the sum \eqref{hessian} has two terms. Since $A$ is invariant under
complex conjugation of the variables, it follows that
$\phi_{1} = - \phi_{2}, \psi_{1} = - \psi_{2}$, hence the two terms
are actually equal. This shows immediately that
$\mu = (\lambda/\pi^2)|_{\am}$ outside $F$. By the following
Lemma neither
$\mu$ nor $\lambda$ has any mass on $F$, so this equality holds everywhere.

\begin{lem}\label{lem3}
If $\Log^{-1}(x) \cap A$ is a finite set for all $x$, then $\mu$ has no mass
on $F$.
\end{lem}

\begin{proof} In Theorem 5 in \cite{PaRu} it is shown that $\mu(E)$
is proportional to
the average number of solutions in $\Log^{-1}(E)$ to the system of equations
\begin{equation}\label{eq2}
f(z_1, z_2) = f(t_1 z_1, t_2 z_2) = 0
\end{equation}
as $(t_1, t_2)$ ranges over the real torus
$\T^2 = \{t \in \C^2; |t_1| = |t_2| = 1\}$. Note that the set of
critical values of the mapping $A \to \R^2 :
(z_1, z_2) \mapsto (|z_1|^2, |z_2|^2)$ is a semialgebraic set. Thus
it is contained in a
real-algebraic curve $\tilde{F}$.

Consider the product space
$\C^2 \times \T^2$ with the two projections $\pi_1$ and $\pi_2$ onto
$\R^2$ and $\T^2$ defined by $\pi_1(z,t) = (|z_1|^2, |z_2|^2)$ and
$\pi_2(z,t) = t$. Let $C = \pi_1^{-1}(\tilde{F}) \cap
\{f(z_1, z_2) = f(t_1 z_1, t_2 z_2) = 0\} \subset \C^2 \times \T^2$.
Since the map $ \pi_1 : C \to \tilde{F}$ has discrete fibers,
it follows that $C$ is a real curve. Hence $\pi_2(C)$ is a null set in $\T^2$.
Since the equation \eqref{eq2} has no solutions in $\Log^{-1}(F)$ for
$t$ outside $\pi_2(C)$, it follows that $\mu(F) = 0$ as required. \end{proof}

\section{Proof of Theorem \ref{main}: conditions 2 and 3 are equivalent.}
%\input{grisha}
%In this section we prove that the conditions 2 and 3 of
%Theorem 1 are equivalent.

\subsection{Implication $2\implies 3$.}
By our assumption $A$ is real up to multiplication
by a constant.
Thus multiplying by a suitable constant we may assume that
$A$ is already defined over $\R$. In this case we may define
the real part $\R A$ as the fixed point set of the
involution of complex conjugation $\conj:(z_1,z_2)\mapsto
(\bar{z_1},\bar{z_2})$ restricted to $A$.

Let $\nu:\tilde{A}\to A$ be the normalization of the curve $A$.
The involution $\conj|_A$ can be lifted to an involution
$\conj_{\tilde{A}}$ on the Riemann surface $\tilde{A}$.
Let $\R\tilde{A}$ be the real part of $\tilde{A}$.
Note that $\nu(\R\tilde{A})\subset\R A$,
but real isolated (singular) points of $\R A$ are not
contained in $\nu(\R\tilde{A})$.

Since $\Log|_A$ is at most 2-1 we can view the map
$\Log\circ\nu:\tilde{A}\to\am$ as a branched double covering.
Let $F\subset\am$ be the branch locus of this covering,
i.e. the set of points whose inverse image under $\Log|_A$
consists of one point.

\begin{lem}
The involution $\conj_{\tilde{A}}$ is the deck transformation
of the branched double covering $\Log\circ\nu$.
\end{lem}
\begin{proof}
The Lemma follows from the fact that $\Log$ maps conjugate points to the
same point, $\Log\circ\conj=\Log$.
\end{proof}

\begin{cor}
$\am=\tilde{A}/\conj_{\tilde{A}}$,
while $F=\Log(\nu(\R\tilde{A}))=\dd\am$.
\end{cor}
\begin{proof}
The curve $\tilde{A}$ is non-singular and therefore
$\tilde{A}/\conj_{\tilde{A}}$ is a smooth surface
with the boundary $\R\tilde{A}$.
\end{proof}

%\begin{cor}
%The non-singular locus of $\R A$ coincides with
%$A\cap\Log^{-1}(\dd\am)$.
%\end{cor}

%Note that $\conj|_A$ reverses the orientation of $A$ and
%$\Log\circ\conj=\Log$.
%Therefore, the set of critical points
%of $\Log|_A$ coincides with $\R A$.
%Thus we get $\R A=A\cap\Log^{-1}(F)$.

%Since $\Log|_A$ is at most 2-1 map of
%degree 0, the curve $F$ is the boundary of $\am$. Indeed,
%$(\Log|_A)^{-1}(F)$ must separate $A$ into two open regions,
%$A_+$, where the local degree of $\Log$ is positive
%(for a choice of an orientation on $\R^2$) and $A_-$,
%where it is negative. Each region $A_\pm$ maps diffeomorphically
%onto $\operatorname{int}(\am)$.
%\mnote{Why does it map onto $\operatorname{int}(\am)$? Couldn't there
%be critical points which are not stable under perturbations?}

%\begin{lem}
%$\Log(\R A)=\dd\am$
%\end{lem}
%\begin{proof}
%Since $\Log\circ\conj=\Log$ and $\Log|_A$ is at most 2-1,
%the inverse image of a point from $\Log(\R A)$ consists of one point.
%Therefore, $\Log(\R A)\subset F=\dd\am$.
%\mnote{$\Log(A)$ should be $\Log(\R A)$ I believe}
%On the other hand, $\Log(\R A)\supset\dd\am$ since
%$\Log^{-1}(\dd \am) \cap A$
%is invariant under the complex conjugation $\conj$.
%\mnote{Do you mean that $\Log^{-1}(\dd \am) \cap A$ is invariant
%under $\conj$?}
%\end{proof}

Thus $\dd\am$ consists of the images of components of $\R\tilde{A}$.
These components are of two types, closed components, called {\em ovals},
and non-compact components.
Accordingly, each oval of $\R A$ which does not contain
singular points corresponds to a hole in $\am$.

Consider first the case when $A$ is a non-singular curve,
so that $\tilde{A}=A$.
Let $l$ be the number of ovals of $\R A$.
Then $\chi(\am)=1-l$, where $\chi$ stands for the homology
Euler characteristic, i.e. the alternated sum of Betti numbers
(we specify that since $\am$ is not compact).
On the other hand, by additivity of Euler characteristic for
compact spaces, $\chi(\bar{A})=2\chi(\am)=2-2l$
(recall that $\bar{A}$ is a compactification of $A$ in
a suitable toric surface, see Section \ref{harnackcurves}).
But $\chi(\bar{A})=2-2g$ and, therefore, $l=g$.

To ensure that $\R A$ has the right number of non-compact
components we recall that $\bar{A}$ intersect the complexification
of $l_k$ in $d_k$ points. Each such intersection corresponds to
a "tentacle" of $\am$ which goes to infinity (see \cite{GKZ}).
Therefore $\R^2\setminus\am$ has $\sum\limits_{k=1}^n d_k$ non-compact
components and each of them must be bounded by a non-compact
component of $\R A$.

To finish the proof in the case when $A$ is non-singular
we need to show that these
$g+\sum\limits_{k=1}^n d_k$ components of $\R A$ are
arranged in $\rtor$ in the Harnack way. This follows
from Lemma 11 of \cite{Mi}. Compactifying with $l_1\cup\dots
\cup l_n$ we obtain that $(\R T_\Delta;\R\bar{A},l_1\cup\dots\cup l_n)$
is a Harnack arrangement.

Now we consider a general case where $A$ might have
singular points.
\begin{lem}
$A$ has no singularities other than real
isolated double points.
\end{lem}
\begin{proof}
We claim that the singular points of $A$ may only
arise as the intersection points of two non-singular
branches of $\tilde{A}$.
Consider the map $\tilde{A}\to A\to \am$.

Over $\am\setminus F$ each of the two branches of $\tilde{A}$
must be non-singular. Indeed, it maps 1-1 to $\am\setminus F$
and, therefore,
the link of each point of this branch is an unknot.

By a similar reason branches of $\tilde{A}$
cannot have singular points over $F$.
Indeed, the links of such points are unknots
since neighborhoods of those points map 2-1 to small half-disks
from $\tilde{A}/\conj_{\tilde{A}}$.

By Lemma 1 of \cite{Mi} the image of each branch of $\tilde{A}$
under $\Log$ has a convex complement. Therefore the images of
branches of $\R\tilde{A}$ cannot intersect (that would produce
points of $\am$ with at least 4 inverse images under $\Log|_A$).

Thus the only singularities of $A$ are intersection points $p$
of a pair of conjugate non-singular imaginary branches.
If these branches are not transverse then they have
a real tangent line $\tau$. The points of $\tau$ close to $p$
will be covered at least twice by each of the two branches of
$\tilde{A}$ which leads to a contradiction.
We conclude that the only singularities of $A$ are $A_1^+$-singularities.
\end{proof}

Now we may replace each $A_1^+$-point with a small oval that
corresponds to its local perturbation and proceed similar to
the case of non-singular curves.

\subsection{Implication $3\implies 2$.}
This implication is contained in the proof of the
main theorem in \cite{Mi}. Indeed, a Harnack curve
is in cyclically maximal position (see Theorem 3 of \cite{Mi}).
By Lemmas 5 and 8 of \cite{Mi} we know that $F=\Log(\R A)=\dd\am$ and
by Lemma 9 $\Log|_{\R A}$ is an embedding.
Therefore the only singularities of $\Log|_A$ are folds and
$\Log|_A$ is at most 2-1.

\end{document}